\documentclass[10pt,reqno]{amsart}
\usepackage{amssymb,amscd,amsbsy}
\renewcommand{\a}{\alpha}

\newcommand{\e}{\varepsilon}

\newcommand{\z}{\zeta}

\renewcommand{\t}{\tau}

\newcommand{\f}{\varphi}
\renewcommand{\o}{\omega}

\renewcommand{\L}{\Lambda}

\renewcommand{\O}{\Omega}

\newcommand{\I}{{\mathcal I}}

\newcommand{\T}{{\Bbb T}}

\newcommand{\R}{{\Bbb R}}

\newcommand{\bS}{{\boldsymbol S}}

\newcommand{\rf}[1]{(\ref{#1})}

\newcommand{\df}{\stackrel{\mathrm{def}}{=}}
\newcommand{\dist}{\operatorname{dist}}

\newcommand{\rank}{\operatorname{rank}}
\newcommand{\const}{\operatorname{const}}

\newcommand{\eeq}{\end{equation}}
\newcommand{\beq}{\begin{equation}}
\newcommand{\bay}{\begin{eqnarray}}
\newcommand{\ba}{\begin{align*}}
\newcommand{\ea}{\end{align*}}
\newcommand{\ey}{\end{eqnarray}}
\newcommand{\bey}{\begin{eqnarray*}}
\newcommand{\eey}{\end{eqnarray*}}

\newcommand{\be}{\infty}

\newcommand{\bl}{\blacksquare}

\newcommand{\Pf}{{\bf Proof. }}

\newcommand{\ov}{\overline}

\newtheorem{thm}{\hspace{\parindent}Theorem}[section]

\begin{document}

\newcommand{\vse}{\vspace{.2in}}

\title{Lipschitz functions of perturbed operators}
\author{Fedor Nazarov and Vladimir Peller}


\maketitle

\footnotesize

\noindent
{\bf Abstract.}

 We prove that if $f$ is a Lipschitz function on $\R$, $A$ and $B$ are self-adjoint operators such that $\rank(A-B)=1$,
then $f(A)-f(B)$ belongs to the weak space $\bS_{1,\be}$, i.e., $s_j(A-B)\le\const(1+j)^{-1}$. We deduce from this result that if $A-B$ 
belongs to the trace class $\bS_1$ and $f$ is Lipschitz, then $f(A)-f(B)\in\bS_\O$, i.e.,
$\sum_{j=0}^ns_j(f(A)-f(B))\le\const\log(2+n)$. We also obtain more general results about the behavior 
of double operator integrals of the form $Q=\iint(f(x)-f(y))(x-y)^{-1}dE_1(x)TdE_2(y)$, where $E_1$ and $E_2$ are spectral measures. We show that 
if $T\in\bS_1$, then $Q\in\bS_\O$ and if $\rank T=1$, then
$Q\in\bS_{1,\be}$. Finally, if $T$ belongs to the Matsaev ideal $\bS_\o$, then $Q$ is a compact operator.

\medskip

\medskip

\noindent
{\bf R\'esum\'e.} 

{\bf Fonctions lipschitziennes d'op\'erateurs perturb\'es.} Nous d\'emontrons que si $f$ est une fonction lipschitzienne, $A$ et $B$   des op\'erateurs autoadjoints tels que $\rank(A-B)=1$, alors $f(A)-f(B)\in\bS_{1,\be}$, c'est-\`a-dire $s_j(A-B)\le\const(1+j)^{-1}$. Si $ A-B $ est dans la classe $\bS_1$ des op\'erateurs \`a trace, nous montrons que 
$f(A)-f(B)\in\bS_\O$, c'est-\`a-dire $\sum_{j=0}^ns_j(f(A)-f(B))\le\const\log(2+n)$.   Plus g\'en\'eralement,  pour une fonction lipschitzienne $f$ et pour des mesures spectrales $E_1$ et $E_2$,  consid\'erons l'int\'egrale double op\'eratorielle 
$Q=\iint(f(x)-f(y))(x-y)^{-1}dE_1(x)TdE_2(y)$. Nous montrons que si $T\in\bS_1$, alors $Q\in\bS_\O$ et si $\rank T=1$, alors $Q\in\bS_{1,\be}$.
Finalement, si $T$ appartient \`a l'id\'eal de Matsaev $\bS_\o$, alors  $Q$ est un op\'erateur compact.

\normalsize

\medskip

\begin{center}
{\bf\large Version fran\c caise abr\'eg\'ee}
\end{center}

\medskip

Dans cette note nous consid\'erons les propri\'et\'es de $f(A)-f(B)$, o\`u $f$ est une fonction lipschitzienne sur la droite r\'eelle $\R$, $A$ et $B$ sont des
op\'erateurs autoadjoints (pas n\'ecessairement born\'es) dont la diff\'erence $A-B$ est ``petite''. Il est bien connu que si $A-B$ appartient \`a l'espace $\bS_1$ des op\'erateurs nucl\'eaires, l'op\'erateur $f(A)-f(B)$  n'appartient pas n\'ecessairement \`a $\bS_1$. 

Nous d\'emontrons que si $A-B\in\bS_1$ et $f$ est une fonction lipschitzienne, alors $f(A)-f(B)$ appartient \`a l'id\'eal $\bS_\O$ d\'efini comme l'ensemble d'op\'erateurs $T$ dont les nombres singuliers $s_j(T)$ satisfont \`a l'in\'egalit\'e
$$
\sum_{j=0}^ns_j(T)\le\const\log(2+n),\quad n\ge0.
$$
Pour d\'emontrer ce r\'esultat nous utilisons la formule de Birman et Solomyak
$$
f(A)-f(B)=\iint\frac{f(x)-f(y)}{x-y}\,dE_A(x)(A-B)\,dE_B(y),
$$
o\`u $E_A$ et $E_B$ sont les mesures spectrales des op\'erateurs $A$ et $B$ (la th\'eorie des int\'egrales doubles op\'eratorielles est d\'evelopp\'ee
dans les travaux \cite{BS1}, \cite{BS2} et \cite{BS3} de Birman et Solomyak). 
Nous \'etablissons un r\'esultat plus g\'en\'eral:  si $f$ est une fonction lipschitzienne, $E_1$ et $E_2$   des mesures spectrales et $T$   un op\'erateur de la classe $\bS_1$, alors
$$
\iint\frac{f(x)-f(y)}{x-y}\,dE_1(x)T\,dE_2(y)\in\bS_\O.
$$

Nous pouvons am\'eliorer les r\'esultats ci-dessus dans le cas $\rank T=1$. En r\'ealit\'e, dans ces cas
$$
\iint\frac{f(x)-f(y)}{x-y}\,dE_1(x)T\,dE_2(y)\in\bS_{1,\be}\df \Big\{T:~\|\T\|_{\bS_{1,\be}} \df\sup_{j\ge0}s_j(T)(1+j)<\be\Big\}.
$$
  Ce fait implique que si $A$ et $B$ sont des op\'erateurs autoadjoints tels que $\rank(A-B)=1$, alors $f(A)-f(B)\in\bS_{1,\be}$.

En utilisant des arguments de dualit\'e on peut montrer que si $T$ appartient \`a l'id\'eal de Matsaev $\bS_\o$, c'est-\`a-dire
$$
\sum_{j\ge0}\frac{s_j(T)}{1+j}<\be,
$$
alors $\iint(f(x)-f(y))(x-y)^{-1}dE_1(x)TdE_2(y)$ est un op\'erateur compact. En particulier, si $A$ et $B$ sont des op\'erateurs autoadjoints tels que
$A-B\in\bS_\o$, alors $f(A)-f(B)$ est un op\'erateur compact. 

Pour \'etablir les r\'esultats ci-dessus nous montrons que si $\mu$ et $\nu$ sont des mesures bor\'eliennes finies sur $\R$, $\f\in L^2(\mu)$, $\psi\in L^2(\nu)$,
$$
k(x,y)=\f(x)\frac{f(x)-f(y)}{x-y}\psi(y),\quad x,\,y\in\R,
$$
et si $\I_k:L^2(\nu)\to L^2(\mu)$ est l'op\'erateur int\'egral d\'efini par $(\I_k g)(x)=\int k(x,y)g(y)\,d\nu(y)$, alors 
$$
\sup_{j\ge0}(1+j)s_j(\I_k)\le\const\|f\|_{\rm Lip}\|f\|_{L^2(\mu)}\|\psi\|_{L^2(\nu)}.
$$

En utilisant des arguments d'interpolation on peut d\'emontrer que si $T$ appartient \`a la classe de Schatten--von Neumann $\bS_p$,
$1\le p<\be$, et $\e>0$, alors 
$$
\iint(f(x)-f(y))(x-y)^{-1}dE_1(x)TdE_2(y)\in\bS_{p+\e}.
$$
En particulier, si $A$ et $B$ sont des op\'erateurs 
autoadjoints tels que $A-B\in\bS_p$, alors $f(A)-f(B)\in\bS_{p+\e}$.

La question de savoir si la condition $T\in\bS_1$ implique que 
$$
\iint(f(x)-f(y))(x-y)^{-1}dE_1(x)TdE_2(y)\in\bS_{1,\be}
$$
est toujours ouverte. Une r\'eponse positive impliquerait que, dans le cas $1<p<\be$,  on a $f(A)-f(B)\in\bS_p$ pour toute paire d'op\'erateurs autoadjoints $A,B$ 
dont la diff\'erence $A-B$ appartient \`a $\bS_p$.

Finalement nous voudrions signaler qu'on peut obtenir des r\'esultats similaires pour les fonctions d'op\'erateurs unitaires et pour les fonctions de contractions.

\medskip

\begin{center}
------------------------------
\end{center}

\setcounter{section}{0}
\section{\bf Introduction}

\medskip

In this note we study the behavior of Lipschitz functions of perturbed operators. It is well known that if $f\in{\rm Lip}$, i.e., $f$ is a Lipschitz function and
$A$ and $B$ are self-adjoint operators with difference in the trace class $\bS_1$, then $f(A)-f(B)$ does not have to belong to $\bS_1$. 
The first example of such $f$, $A$, and $B$ was constructed in \cite{F}. Later in \cite{Pe1} a necessary condition on $f$ was found ($f$ must be locally in the Besov space $B_1^1$) under which 
the condition $f(A)-f(B)\in\bS_1$ implies that $f(A)-f(B)\in\bS_1$. That necessary condition also implies that the condition $f\in{\rm Lip}$ is not sufficient.

On the other hand, Birman and Solomyak showed in \cite{BS3} that if $A-B$ belongs to the Hilbert--Schmidt class $\bS_2$, then
$f(A)-f(B)\in\bS_2$ and $\|f(A)-f(B)\|_{\bS_2}\le\|f\|_{\rm Lip}\|A-B\|_{\bS_2}$, where 
$\|f\|_{\rm Lip}\df\sup_{x\ne y}|f(x)-f(y)|\cdot|x-y|^{-1}$.
Moreover, it was shown in \cite{BS3} that in this case $f(A)-f(B)$ can be expressed in terms of the following double operator integral
\bay
\label{BSf}
f(A)-f(B)=\iint\frac{f(x)-f(y)}{x-y}\,dE_A(x)(A-B)\,dE_B(y).
\ey
where $E_A$ and $E_B$ are the spectral measures of $A$ and $B$. We refer the reader to \cite{BS1}, \cite{BS2}, and \cite{BS3} for the beautiful theory of double operator integrals.
Note that the divided difference $(f(x)-f(y))/(x-y)$ is not defined on the diagonal. Throughout this note we assume that it is zero on the diagonal.

In this note we study properties of the operators $f(A)-f(B)$ for (not necessarily bounded) self-adjoint operators $A$ and $B$ such that $A-B$ has rank one or $A-B\in\bS_1$. Actually, we consider more general operators of the form 
\bay
\label{DOI}
\I_{E_1,E_2}(f,T)\df
\iint\frac{f(x)-f(y)}{x-y}\,dE_1(x)T\,dE_2(y),
\ey
where $E_1$ and $E_2$ are Borel spectral measures on $\R$ and $\rank T=1$ or $T\in\bS_1$. Duality arguments also allow us to study
double operator integrals \rf{DOI} in the case when $T$ belongs to the {\it Matsaev ideal} $\bS_\o$.

Recall the definitions of the following operator ideals:
$$
\bS_{1,\be}\df\Big\{T:~\|\T\|_{\bS_{1,\be}} \df\sup_{j\ge0}s_j(T)(1+j)<\be\Big\},
$$
$$
\bS_\O\df\Big\{T:~\|T\|_{\bS_\O}\df\big(\log(2+n)\big)^{-1}\sum_{j=0}^ns_j(T)<\be\Big\},
$$
and
$$
\bS_\o\df\Big\{T:~\|T\|_{\bS_\o}\df\sum_{j=0}^\be \frac{s_j(T)}{1+j}<\be\Big\}.
$$
It is well known that $\bS_{1,\be}$ is not a Banach space and its Banach hull coincides with $\bS_\O$. Also recall that the dual space to $\bS_\o$ can be identified in a natural way with $\bS_\O$.

Note that the recent paper \cite{AP} contains results on properties of $f(A)-f(B)$  for $f$ in the H\"older class $\L_\a$, $0<\a<1$, and
self-adjoint operators $A$ and $B$ with $A-B$ in Schatten--von Neuman classes $S_p$.

\medskip

\section{\bf Main results}

\medskip

\begin{thm}
\label{r1}
Let $f\in{\rm Lip}$ and let $E_1$ and $E_2$ be Borel spectral measures on $\R$. If $\rank T=1$, then
$\I_{E_1,E_2}(f,T)\in\bS_{1,\be}$ and
$$
\|\I_{E_1,E_2}(f,T)\|_{\bS_{1,\be}}\le\const\|f\|_{\rm Lip}\|T\|.
$$
\end{thm}

Theorem \ref{r1} immediately implies the following result.

\begin{thm}
\label{S1}
Let $f\in{\rm Lip}$ and let $E_1$ and $E_2$ be Borel spectral measures on $\R$. If $T\in\bS_1$, then
$\I_{E_1,E_2}(f,T)\in\bS_\O$ and
$$
\|\I_{E_1,E_2}(f,T)\|_{\bS_\O}\le\const\|f\|_{\rm Lip}\|T\|_{\bS_1}.
$$
\end{thm}

By duality, we obtain the following theorem.

\begin{thm}
\label{Ma}
Let $f\in{\rm Lip}$, and let $E_1$ and $E_2$ be Borel spectral measures on $\R$. Then
the transformer $T\mapsto\I_{E_1,E_2}(f,T)$ defined on $\bS_2$ extends to a bounded linear operator from 
$\bS_\o$ to the ideal of all compact operator and
$$
\|\I_{E_1,E_2}(f,T)\|\le\const\|f\|_{\rm Lip}\|T\|_{\bS_\o}.
$$
\end{thm}

Using interpolation arguments, we can easily obtain from Theorem \ref{S1} the following fact.

\begin{thm}
\label{ep}
Let $f\in{\rm Lip}$, and let $E_1$ and $E_2$ be Borel spectral measures on $\R$. Suppose that $1\le p<\be$ and $\e>0$.
If $T\in\bS_p$, then $\I_{E_1,E_2}(f,T)\in\bS_{p+\e}$.
\end{thm}

Birman--Solomyak formula \rf{BSf} allows us to deduce straightforwardly from Theorems \ref{r1}, \ref{S1}, and \ref{Ma}
the following theorem.

\begin{thm}
\label{one}
Let $A$ and $B$ be self-adjoint operators on Hilbert space and let $f\in{\rm Lip}$. We have

{\em (i)} if $\rank(A-B)=1$, then $f(A)-f(B)\in\bS_{1,\be}$ and $\|f(A)-f(B)\|_{\bS_{1,\be}}\le\const\|f\|_{\rm Lip}\|A-B\|$;

{\em(ii)} if $A-B\in\bS_1$, then $f(A)-f(B)\in\bS_\O$ and $\|f(A)-f(B)\|_{\bS_\O}\le\const\|f\|_{\rm Lip}\|A-B\|_{\bS_1}$;

{\em(iii)} if $A-B\in\bS_\o$, then $f(A)-f(B)$ is compact and $\|f(A)-f(B)\|\le\const\|f\|_{\rm Lip}\|A-B\|_{\bS\o}$;

{\em(iv)} if $1\le p<\be$, $\e>0$, and $A-B\in\bS_p$, then $f(A)-f(B)\in\bS_{p+\e}$.
\end{thm}

It is still unknown whether the assumption $T\in\bS_1$ implies that $\I_{E_1,E_2}(f,T)\in\bS_{1,\be}$. If this is true, then the condition $A-B\in\bS_p$ would imply that $f(A)-f(B)\in\bS_p$ for $1<p<\be$.

To prove Theorem \ref{r1}, we obtain a weak type estimate for Schur multipliers.

For a kernel function $k\in L^2(\mu\times\nu)$, we define the integral operator $\I_k:L^2(\nu)\to L^2(\mu)$ by
$$
(\I_k g)(x)=\int k(x,y)g(y)\,d\nu(y),\quad g\in L^2(\nu).
$$

As in the case of transformers from $\bS_1$ to $\bS_1$ (see \cite{BS3}), Theorem \ref{r1} reduces to the following fact.

\begin{thm}
\label{Sm}
Let $\mu$ and $\nu$ be finite Borel measures on $\R$, $\f\in L^2(\mu)$, $\psi\in L^2(\nu)$.
Suppose that $f\in{\rm Lip}$ and the kernel function $k$ is defined by
$$
k(x,y)=\f(x)\frac{f(x)-f(y)}{x-y}\psi(y),\quad x,\,y\in\R.
$$
Then the integral operator $\I_k: L^2(\nu)\to L^2(\mu)$ with kernel function $k$ belongs to $\bS_{1,\be}$ and
$$
\|\I_k\|_{\bS_{1,\be}}\le\const\|f\|_{\rm Lip}\|\f\|_{L^2(\mu)}\|\psi\|_{L^2(\nu)}.
$$
\end{thm}

\Pf Without loss of generality we may assume that
$\|\f\|_{L^2(\mu)}=\|\psi\|_{L^2(\nu)}=1$ and $\|f\|_{\rm Lip}=1$.
Let us fix a positive integer $n$. 

Given $N>0$, we denote by $P_N$ multiplication by 
the characteristic function of $[-N,N]$ (we use the same notation for multiplication on $L^2(\mu)$ and on $L^2(\nu)$).
Then for sufficiently large values of $N$,
\bay
\label{N}
\|\I_k-P_N\I_kP_N\|_{\bS_2}<\frac1{n^{1/2}}.
\ey
Clearly, $P_N\I_kP_N$ is the integral operator with kernel function $k_N$, 
$k_N(x,y)=\chi_N(x)k(x,y)\chi_N(y)$,
where $\chi_N=\chi_{[-N,N]}$ is the characteristic function of $[-N,N]$. We fix $N>0$, for which \rf{N} holds.

Consider now the points $x_j$, $1\le j\le r$, and $y_j$, $1\le j\le s$, at which $\mu$ and $\nu$ have point masses and
\bay
\label{mass}
|\f(x_j)|^2\mu\{x_j\}\ge\frac1n,\quad1\le j\le r,\quad\mbox{and}\quad|\psi(y_j)|^2\nu\{y_j\}\ge\frac1n,\quad1\le j\le s.
\ey
Clearly, $r\le n$ and $s\le n$.
We define now the kernel function $k_\sharp$ by
$$
k_\sharp(x,y)=u(x)k_N(x,y)v(y),\quad x,\,y\in\R,
$$
where
$$
u(x)\df1-\chi_{\{x_1,\cdots,x_r\}}(x)\quad\mbox{and}\quad v(y)\df1-\chi_{\{y_1,\cdots,y_s\}}(y).
$$
Obviously, the integral operators $\I_{k_N}$ and $\I_{k_\sharp}$ coincide on a subspace of codimension at most $r+s\le2n$.

We can split now the interval $[-N,N]$ into no more than $n$ subintervals $I$, $I\in{\frak I}$, such that
$$
\int_I|\f(x)|^2u(x)\,d\mu(x)+\int_I|\psi(y)|^2v(y)\,d\nu(y)\le\frac4n,\quad I\in{\frak I}.
$$
This is certainly possible because of \rf{mass}.

We have
$\I_{k_\sharp}=\I^{(1)}+\I^{(2)}+\I^{(3)}$,
where
$$
\big(\I^{(1)}g\big)(x)=\int\limits_\R\left(\sum_{I\in{\frak I}}\chi_I(x)k_\sharp(x,y)\chi_I(y)\right)g(y)\,d\nu(y),
$$
$$
\big(\I^{(2)}g\big)(x)=\int\limits_\R\left(\sum_{I,J\in{\frak I},\,I\ne J,\,|I|\ge|J|}\chi_I(x)k_\sharp(x,y)\chi_I(y)\right)g(y)\,d\nu(y),
$$
and
$$
\big(\I^{(3)}g\big)(x)=\int\limits_\R\left(\sum_{I,J\in{\frak I},\,|I|<|J|}\chi_I(x)k_\sharp(x,y)\chi_I(y)\right)g(y)\,d\nu(y)
$$
(we denote by $|I|$ the length of $I$).
It is easy to see that
$\big\|\I^{(1)}\big\|_{\bS_2}\le4n^{-1/2}$.
Let us estimate $\I^{(2)}$. The integral operator $\I^{(3)}$ can be estimated in the same way.

Suppose that $I,\,J\in{\frak I}$, $I\ne J$, and $|I|\ge|J|$. For $x\in I$ and $y\in J$, we have
$$
\frac1{x-y}=\frac1{x-c(J)}+\frac{y-c(J)}{x-c(J)}\cdot\frac1{x-y},
$$
where $c(J)$ denotes the center of $J$.

Suppose that $g\perp\ov{\psi}\chi_J$ and $g\perp\ov{\psi}\bar f\chi_J$. Then
$\I_2g=\I_{k_\flat}g$,
where
$$
k_\flat(x,y)
=\sum_{I,J\in{\frak I},\,I\ne J,\,|I|\ge|J|}u(x)\f(x)a_{IJ}(x,y)\psi(y)v(y)
$$
and
$$
a_{IJ}(x,y)=\chi_I(x)\frac{y-c(J)}{x-c(J)}\cdot\frac{f(x)-f(y)}{x-y}\chi_J(y).
$$
Thus $\I^{(2)}$ and $\I_{k_\flat}$ coincide on a subspace of codimension at most $2n$.

To estimate the Hilbert--Schmidt norm of $\I_{k_\flat}$, we observe that
$$
|a_{IJ}(x,y)|\le\frac{|J|}{\big(|J|+\dist(I,J)\big)},\quad x\in I,~y\in J.
$$
Thus
\begin{align*}
\big\|\I_{k_\flat}\big\|^2_{\bS_2}&\le\sum_{I,J\in{\frak I},\,I\ne J,\,|I|\ge|J|}\left(\int_I|\f|^2u\,d\mu\right)\left(\int_J|\psi|^2v\,d\nu\right)\|a_{IJ}\|^2_{L^\be}\\[.2cm]
&\le\frac4{n^2}\sum_{I,J\in{\frak I},\,I\ne J,\,|I|\ge|J|}\frac{|J|^2}{\big(|J|+\dist(I,J)\big)^2}.
\end{align*}

Let us observe that for a fixed $J\in{\frak I}$,
\bay
\label{es}
\sum_{I\in{\frak I},\,I\ne J,\,|I|\ge|J|}\frac{|J|^2}{\big(|J|+\dist(I,J)\big)^2}\le\const.
\ey
Indeed, we can enumerate the intervals $I\in{\frak I}$ satisfying $I\ne J$ and $|I|\ge|J|$ so that the resulting intervals $I_k$
satisfy 
$\dist(I_k,J)\le\dist(I_{k+1},J)$.
Since the intervals $I_k$ are disjoint, we have
$$
\dist(I_k,J)\ge\frac{k-3}2|J|.
$$
This easily implies \rf{es}.
It follows that
$$
\|\I_{k_\flat}\|^2_{\bS_2}\le C\frac4{n^2}\cdot n=\frac{4C}n.
$$
Similarly, $\I^{(3)}$ coincides on a subspace of codimension at most $2n$ with an operator whose Hilbert--Schmidt norm is at most
$2\left(C/n\right)^{1/2}$.

If we summarize the above, we see that $\I_k$ coincides on a subspace of codimension at most $6n$ with an operator whose Hilbert--Schmidt norm
is at most $Kn^{-1/2}$, where $K$ is a constant.
Hence, on a subspace of codimension at most $7n$ the operator $\I_k$ coincides with an operator whose norm is at most
$K/n$, i.e.,
$$
s_{7n}(\I_k)\le\frac Kn,\quad n\ge1,\quad\bl
$$

 Note that in the case of operators on the space $L^2(\T)$ with respect to Lebesgue measure on the unit circle $\T$,  
the following related fact was obtained in \cite{Pe} (see also \cite{Pe3}): if the derivative of $f$ belongs to the Hardy class $H^1$,  $\f$ and $\psi$ belong to $L^\be(\T)$, and the kernel function $k$ is defined  by
$$
k(\z,\t)=\f(\z)\frac {f(\z)-f(\t)}{\z-\t}\psi(\t),\quad\z,\,\t\in\T,
$$
then the integral operator $\I_k$ on $L^2(\T)$ belongs to $\bS_{1,2}$, i.e., $\sum_{j\ge0}(s_j(\I_k))^2(1+j)<\be$.

To conclude the article, we note that similar results can be obtained for functions of unitary operators and for functions of contractions.

{\bf Remark.} After this article had been written we have been informed by D. Potapov and F. Sukochev that they had proved the following result:
if $f$ is a Lipschitz function, $1<p<\be$, and $A$ and $B$ are self-adjoint operators such that $A-B\in\bS_p$, then
$f(A)-f(B)\in\bS_p$.

\footnotesize

\noindent
\begin{tabular}{p{10cm}p{15cm}}
F.L. Nazarov & V.V. Peller \\
Department of Mathematics & Department of Mathematics \\
University of Wisconsin  & Michigan State University \\
Madison, WI 53706 & East Lansing, Michigan 48824\\
USA&USA
\end{tabular}

\end{document}